\documentclass[11pt,reqno]{amsart}
\usepackage{graphicx}
\pdfoutput=1
\usepackage{amsfonts,amsmath,amssymb}
\usepackage{hyperref}
\usepackage[linesnumbered,ruled,vlined,norelsize]{algorithm2e}
\usepackage{amsthm}
\usepackage{amsmath}
\usepackage{amssymb}
\let\amslrcorner\lrcorner
\usepackage[utf8]{inputenc}
\usepackage{geometry}
\usepackage{hyperref}
\usepackage{mathtools}
\usepackage{calrsfs}
\usepackage{epigraph}
\usepackage{tikz-cd}
\usepackage{mathabx}
\usepackage{tcolorbox}
\usepackage{subfiles}
\usepackage{commath}
\usepackage{mathtools}
\let\lrcorner\amslrcorner
\usepackage{musicography}
\usepackage{scalerel,stackengine}

\usepackage[backend=biber]{biblatex}
\usepackage{biblatex}
\theoremstyle{plain}
\newtheorem{Thm}{Theorem}[section]
\newtheorem{Prop}{Proposition}[section]

\newtheorem{question}{Question}[section]

\theoremstyle{definition}
\newtheorem{Def}{Definition}[section]
\newtheorem{Exa}{Example}[section]

\theoremstyle{remark}
\newtheorem{Rem}{Remark}[section]
\numberwithin{equation}{subsection}

\numberwithin{equation}{section}

\makeatletter
\DeclareFontFamily{OMX}{MnSymbolE}{}
\DeclareSymbolFont{MnLargeSymbols}{OMX}{MnSymbolE}{m}{n}
\SetSymbolFont{MnLargeSymbols}{bold}{OMX}{MnSymbolE}{b}{n}
\DeclareFontShape{OMX}{MnSymbolE}{m}{n}{
    <-6>  MnSymbolE5
   <6-7>  MnSymbolE6
   <7-8>  MnSymbolE7
   <8-9>  MnSymbolE8
   <9-10> MnSymbolE9
  <10-12> MnSymbolE10
  <12->   MnSymbolE12
}{}
\DeclareFontShape{OMX}{MnSymbolE}{b}{n}{
    <-6>  MnSymbolE-Bold5
   <6-7>  MnSymbolE-Bold6
   <7-8>  MnSymbolE-Bold7
   <8-9>  MnSymbolE-Bold8
   <9-10> MnSymbolE-Bold9
  <10-12> MnSymbolE-Bold10
  <12->   MnSymbolE-Bold12
}{}

\let\llangle\@undefined
\let\rrangle\@undefined
\DeclareMathDelimiter{\llangle}{\mathopen}%
                     {MnLargeSymbols}{'164}{MnLargeSymbols}{'164}
\DeclareMathDelimiter{\rrangle}{\mathclose}%
                     {MnLargeSymbols}{'171}{MnLargeSymbols}{'171}
\makeatother

\renewcommand{\rm}{\normalshape}

\newcommand{\ds}{\displaystyle}
\newcommand{\C}{{\mathbb C}}
\newcommand{\R}{{\mathbb R}}

\newcommand{\VV}{{\mathcal{V}}}

\newcommand{\dd}{{\mathrm{d}}}

\newcommand{\vspp}{\vspace*{5pt}}

\newcommand{\IM}{{\mathsf{Im}\,}}
\newcommand{\RE}{{\mathsf{Re}\,}}

\newcommand{\til}{\widetilde}
\DeclareMathOperator{\spn}{span}
\newcommand{\brr}{\overline}

\DeclareMathOperator{\aut}{aut}

\DeclareMathOperator{\hol}{hol}

\def \im{\text{\rm Im }}

\newcommand{\CC}[1]{\mathbb{C}^{#1}}

\begin{document}

\author{Ilya Kossovskiy}
\address{Department of Mathematics $\&$  International Center of Mathematics,  Southern University of Science and Technology, Shenzhen, China}

\address{Department of Mathematics, Masaryk University, Brno, Czechia}

\address{Institute of Discrete Mathematics and Geometry, Vienna University of Technology, Vienna, Austria}

\email{ilyakos@sustech.edu.cn}

\author{Vinícius Novelli}
\address{Department of Mathematics, Instituto de Ciências Matemáticas e de Computação, Universidade de São Paulo, São Carlos, Brazil}
\email{viniciusnovelli@usp.br}

\title{On the regularity of nondegenerate hypo-analytic structures of hypersurface type}\footnote{MSC 35F05, 32V05, 32V25}

\date{January 28, 2025.
        $^*$Corresponding author: viniciusnovelli@usp.br}


\maketitle

\begin{abstract} For a smooth, non-degenerate locally integrable structure of hypersurface type on a manifold $M$, we provide necessary and sufficient conditions for it to be equivalent, near a point, to a real-analytic locally integrable structure ({\em the analytic regularizability}), generalizing a recent result of Zaitsev and the first author \cite{koss_zaitsev}. 

First, we discover, in our setting, a (previously unknown) invariant CR submanifold $\Sigma$ in $M$ of  hypersurface type, which we call {\em the central submanifold}. 
We prove that the analytic regularizability of $M$ is equivalent to that of the associated CR manifold $\Sigma$. Furthermore, as a byproduct of our construction, we show that the central manifold construction reduces the whole (smooth or analytic) equivalence problem  for nondegenerate structures with the Levi positivity condition to that of the associated central manifolds, i.e. to CR geometry.

Second, we make use of a classical construction due to Marson \cite{Marson} and show that sufficient for the analytic regularizability of $M$ is the analytic regularizability of the CR manifold $\tilde M$ associated with $M$ in the sense of Marson.

We show applications of both regularizability conditions to classes of locally integrable structures.
\end{abstract}

\tableofcontents
\section{Introduction}

This work is concerned with {\em locally integrable structures} of hypersurface type. These are locally integrable systems of complex vector fields on a smooth manifold for which the rank of the characteristic set does not exceed one. Within this class are, in particular, CR structures on real hypersurfaces in a complex manifold and (locally integrable) Mizohata structures.  The theory of locally integrable structures has a long and rich history (see for instance the classical sources Berhanu-Cordaro-Hounie \cite{BCH_book}, Treves \cite{T1}, Cordaro-Treves  \cite{CordTreves}, Baouendi-Chang-Treves \cite{BaouendiChangTreves}). Such structures can be of several regularity levels ($C^k$, smooth, Gevrey, real-analytic), and as a rule, the higher the regularity is, the more techniques are available to understand them. For example, in the real-analytic case we have such powerful techniques  as the {\em external complexification}, {\em Segre varieties} (see, for example, Baouendi-Ebenfelt-Rothschild \cite{BER}), and {\em normal forms} (e.g. Chern-Moser \cite{chern}, Kolar-Kossovskiy-Zaitsev \cite{kkz}).  The latter techniques are available solely in the analytic case, which makes the existence of coordinates where a structure becomes real-analytic a crucial property. 

This inspires the following widely open question, which is our main concern in this paper:

\begin{question}\label{question1} Let $(M,\mathcal{V},p)$ be a germ of a smooth locally integrable structure (of hypersurface type) near a point $p\in M$. Find necessary and sufficient conditions for $\mathcal{V}$ to be equivalent to a germ of a \textit{real-analytic} locally integrable structure.
\end{question}

In the situation when the regularity property from \autoref{question1} holds, we say that $\mathcal{V}$ is \textit{analytically regularizable} near $p\in M$. 

We shall point out  some subtlety with the notion of equivalence between integrable structures. In this paper, we consider it in the weak sense: an {\em equivalence between structures} is a germ of a diffeomorphism (of certain regularity) for which the derivative preserves the bundle. There is also the notion of \textit{hypo-analytic equivalence} between these structures (once a hypo-analytic structure is fixed, defined by a system of first integrals). We do not pursue this problem in this paper. These notions will be recalled in Section \ref{preliminaries}.

For some classes of structures, Question \ref{question1} is well understood. For instance, if the structure $\mathcal{V}$ is \textit{elliptic} at $p$ (meaning that the characteristic set at this point is trivial), then, by a theorem of Nirenberg \cite{Nirenberg} , one can find coordinates $(x_1,\ldots,x_m,y_1,\ldots,y_m,t_1,\ldots,t_{n'})$ near $p$ such that a complete system of first integrals for $\mathcal{V}$ is given by $z_j=x_j+iy_j$, $j=1,\ldots,m$. In particular, all such structures are analytically regularizable. Another case where the result is known is that of Mizohata structures (when the locally integrable structure has corank $1$), when an argument using Morse's lemma yields coordinates $(x,t_1,\ldots,t_n)$ such that the first integral is of the form $x+i\mathcal{L}(t) $, where $\mathcal{L}$  is a non-degenerate quadratic form in $t$ (see Theorem VII.1.1 in \cite{T1}). In particular, they are also analytically regularizable. 

To the best of our knowledge, for other classes of structures, the regularizability problem as stated in \ref{question1} was widely open until the recent work \cite{koss_zaitsev} by Zaitsev and the first author. They considered the case of smooth strictly pseudoconvex hypersurfaces in $\C^n$. They introduce a property (called \textit{condition} (E)), which we review in detail in Section 2 and mention now that this property is a holomorphic extension property for certain smooth function associated with a real hypersurface $M$ {\em (the $\Phi$-function}). The main result of \cite{koss_zaitsev} is at a glance the following.

\begin{Thm}[Kossovskiy, Zaitsev \cite{koss_zaitsev}]\label{koss_zaitsev_thm} Let $M\subset \C^n$ be a smooth, strictly pseudoconvex hypersurface, and let $p\in M$. Then $M$ is CR equivalent near $p$ to a real-analytic hypersurface if and only if $M$ satisfies condition (E).
\end{Thm}

The goal of this paper is to extend this result to a larger class of locally integrable structures (of hypersurface type). One can view, through the \textit{local fine embedding} (see Section \ref{preliminaries}), such a structure as a smooth family of Levi non-degenerate hypersurfaces, and the issue is to find a smooth family of normalizing transformations satisfying the additional equations coming from the extra parameters.

Our technique yields two Levi non-degenerate hypersurfaces associated to $\mathcal{V}$. One (which is intrinsic) consists of the set of points for which the structure is not elliptic. This submanifold, $\Sigma \subset M$, which we refer to as the \textit{central submanifold of} $\mathcal{V}$, inherits a structure of the ambient space that of a Levi non-degenerate hypersurface, and yields a necessary condition for regularizability of $\mathcal{V}$. We state the first main result, and the details can be found in Section \ref{necessary_section}.
\begin{Thm}\label{main} Let $(M,\mathcal{V},p)$ be a germ of a nondegenerate locally integrable structure of hypersurface type. Assume the Levi form at $p$ is definite. Then,  $\mathcal{V}$ is analitically regularizable if and only if the central CR hypersurface $\Sigma$ satisfies condition (E) at $p$.

\end{Thm}

On the other hand, one can fix a set of first integrals (induced by a particular system of coordinates) and construct, following the idea of Marson in \cite{Marson}, an associated (external) hypersurface, which in the positive-definite Levi form case is a strictly pseudoconvex CR hypersurface. Our second main result is that regularizability of this external hypersurface entails regularizability of $\mathcal{V}$. The details can be found in Section 5.

\begin{Thm}\label{sufficient_thm} Let $(M,\mathcal{V},p)$ be a germ of a nondegenerate locally integrable structure of hypersurface type. Assume that the Levi form at $p$ is definite. Fix a system of coordinates given by the fine local embedding near $p$ and let $M^{\bullet}$ be the external hypersurface associated to $M$ at $p$ in these coordinates. Then, if $M^{\bullet}$ satisfies condition (E), then $\mathcal{V}$ is analytically regularizable at $p$.
\end{Thm}

Finally, we outline that, as a bi-product of our considerations, we show the following remarkable property: {\em if a given nondegenerate structure with a definite Levi form is rigid, then it is analytically regularizable if and only if its central manifold is actually analytic}. This result is new even in the case of CR structures. Details can be found in Section 6.

\begin{Rem}
In an upcoming paper of Zaitsev and the first author \cite{KZnew}, Condition (E) above is replaced by a series of bounds for certain derivatives of the defining function of the hypersurface. For each given hypersurface, these bounds can be directly verified, which further emphasizes the effectiveness of the regularizability conditions for locally integrable structures introduced in the paper. 
\end{Rem}
\begin{Rem}
In fact, the construction of the central manifold associated with a locally integrable structure in our paper gives consequence going much beyond the regularizability problem. Namely (see Theorem \autoref{equivalence} below), it reduces the {\em entire} (smooth or analytic) equivalence problem for nondegenerate structures with the Levi positivity  to that for the associated central manifolds. The latter makes then possible the application of the whole machinery of CR geometry (e.g. normal forms \cite{chern,kkz} and methods of {\em Parabolic Geometries} \cite{CS}) for the classification of  the involutive structures under discussion. 
\end{Rem}

\begin{center}
\bf Acknowledgements

\end{center}

\medskip

The first author was supported by the Sustech University internal Grant Y01286147, NSFC grant K24281006, Austrian Science Fund Grant P34369, and GACR Grant GC22-15012. The second author was supported by FAPESP, grant 2023/17607-7.

The authors would like to thank Paulo Cordaro for numerous fruitful discussions during the preparation of the paper.

\section{Preliminaries}\label{preliminaries}

\subsection{Nondegenerate locally integrable structures.} A \textit{smooth locally integrable structure} is the datum of a smooth $N$-dimensional manifold $M$ and a smooth complex vector subbundle $\mathcal{V}\subset \C T M$ (of complex rank $n$) such that the Lie bracket of two local sections is still a section, and such that the orthogonal $T':=\mathcal{V}^{\perp} \subset \C T^{\ast}M$ (for the duality between vector fields and one-forms) is locally generated by exact one-forms. A distribution $u$ defined on an open set of $M$ is a \textit{solution} of $\mathcal{V}$ if $\dd u$ is a local section of $T'$. A family of complex-valued functions $Z_1,\ldots,Z_m$, where $m=N-n$, is a \textit{complete set of first-integrals} for $\mathcal{V}$ over an open set $U\subset M$ if $\dd Z_1,\ldots,\dd Z_m$ forms a basis of $T'$ over $U$.

Given a point $p\in M$, the \textit{characteristic set} of $\mathcal{V}$ at $p$ is the space $T' \cap T^\ast_pM =: T^{\circ}_p$ of the real covectors of $T'$. We say the structure is \textit{of hypersurface type} if 
\[
\dim_{\R} T^{\circ}_p \leq 1,\,\,\,\,\,p\in M.
\]
It is important to note that, in general, the association $$M\ni p \mapsto T^{\circ}_p$$ does \textit{not} define a vector bundle over $M$. This is true in the case of CR structures, but we shall encounter structures (of Mizohata type) for which the dimension of $T^{\circ}_p$ is not locally constant (it is, however, an upper semicontinuous function of $p$).

\vspp

The results of this paper are local, so we will always reason in a neighborhood of a distinguished point, which we henceforth refer to as the origin and denote by $0$. When $\dim T^{\circ}_0 = 0$, the structure is {\em elliptic} in a neighborhood of $0$, which is not an interesting case for the problems discussed here. We will therefore work under the hypothesis that
\[
\dim_{\R} T^{\circ}_0 = 1.
\]

We can select a solution $w$ in a neighborhood of $0$ such that $\dd w$ spans $T^{\circ}_0$ at the origin. Select $m-1=:\nu$ solutions $z_1,\ldots,z_\nu$ in some open neighborhood of $0$ such that $\dd z_1 \wedge \ldots \wedge \dd z_\nu \wedge \dd w \not = 0$ at $0$. 

It is a matter of linear algebra (see \cite{T1} or \cite{BCH_book}) to see that this is equivalent to 
\[
\dd z_1 \wedge \ldots \wedge \dd z_\nu \wedge \dd \brr{z_1} \wedge \ldots \wedge \dd \brr{z_\nu} \wedge \dd (\RE w) \not= 0,\,\,\,\,\dd (\IM w)=0
\]
at the origin.

We can conveniently take $x_i = \RE z_i$, $y_j = \IM z_j$ and $s=\RE w$ as coordinates (for $i,j\in\{1,\ldots,\nu\}$) and adjoin to them coordinates $t_1,\ldots,t_{n'}$ ($n':=n-\nu$), obtaining a complete system of coordinates for $M$ near $0$. Adding appropriate constants allows us to assume that this system is centered at $0$, and $\IM w(0)=0$. We conclude that the differentials of the functions 
\begin{equation}\label{local_coord}
z_j = x_j+iy_j,\,\,\,\,\,j=1,\ldots,\nu,\,\,\,\,\,\,w=s+i\phi(z,\brr{z},s,t)
\end{equation}
span the bundle $T'$ in a neighborhood of $0$, and moreover we have $\phi(0)=0$ and $\dd \phi(0)=0$. The characteristic set is spanned at the origin by $\dd s\big|_{0}$. We can then construct a frame for $\mathcal{V}$ near $0$ in the following way:
\begin{equation}\label{vector_fields}
L_j = \frac{\partial}{\partial \brr{z_j}} - i \frac{\phi_{\brr{z_j}}}{1+i\phi_s}\frac{\partial}{\partial s},\,\,\,j=1,\ldots,\nu,\,\,\,\,\,\,\,L_{\nu+l}=\frac{\partial}{\partial t_l}-i\frac{\phi_{t_l}}{1+i\phi_s}\frac{\partial}{\partial s},\,\,l=1,\ldots,n'.
\end{equation}

One very important invariant of $\mathcal{V}$ is the \textit{Levi form}, which we discuss next. 

\begin{Def} Given a point $p\in M$ and a characteristic covector $\omega \in T^{\circ}_p\setminus \{0\}$, consider
\begin{align*}
\mathcal{B}_{\omega,p}:\mathcal{V}_p \times \mathcal{V}_p &\rightarrow \C \\
(v,w)&\mapsto \frac{1}{2i}\omega\left([L_1\big|_{p},\overline{L_2\big|_{p}}] \right),
\end{align*}
where $L_1,L_2$ are local sections of $\mathcal{V}$ such that $L_1\big|_{p}=v$, $L_2\big|_{p}=w$. The value $\mathcal{B}_{\omega,p}(v,w)$ does not depend on the sections $L_1,L_2$ chosen, and defines a hermitian form on $\mathcal{V}_p$. The associated quadratic form 
\[
\mathcal{L}_{\omega, p}(v):=\mathcal{B}_{\omega,p}(v,v),\,\,v\in \mathcal{V}_p,
\]  
is called the \textit{Levi form} of $\mathcal{V}$ at $(\omega,p)$.
\end{Def}
In the case of structures of hypersurface type, since $\dim_{\R}T^{\circ}_0=1$, it is customary to suppress the dependence on the characteristic covector, and speak only of the Levi form of the structure at $0$. In the coordinates described in \ref{local_coord}, the Levi form (at the characteristic vector $\dd s\big|_{0}$) can be represented by the following $n\times n$ matrix:
\begin{equation}\label{Levi_form}
\mathcal{L}_p = \begin{pmatrix}
\left(\frac{\partial^{2}\phi}{\partial z_j \partial \brr{z_k}}(0) \right)_{\nu\times \nu} & \left(\frac{\partial^{2}\phi}{\partial z_j \partial t_k}(0) \right)_{\nu \times n'} \\
\left(\frac{\partial^2 \phi}{\partial t_j \partial \brr{z_k}}(0) \right)_{n'\times \nu} & \left(\frac{\partial^2 \phi}{\partial t_j \partial t_k}(0) \right)_{n'\times n'}
\end{pmatrix}.
\end{equation}

The structures that we shall work with in this paper are	 defined, then, as follows.
\begin{Def} A \textit{(Levi-)nondegenerate structure of hypersurface type} $\mathcal{V}$ at $0\in M$ is a locally integrable structure of hypersurface type such that the Levi form is nondegenerate at $0$. A nondegenerate structure of hypersurface type is one which is of this form for all $p\in M$.
\end{Def}

The two main examples of structures of this kind are the following:
\begin{Exa}[Levi nondegenerate hypersurfaces] If $n'=0$, i.e., there are no variables $t_j$, then the structure is a Levi nondegenerate real hypersurface $M \subset \C^{n+1}$, given in a neighborhood of zero by $\IM w = \phi(z,\brr{z},\RE w)$, where coordinates in $\C^{n+1}$ are $(z,w)\in \C^n \times \C$.
\end{Exa}
\begin{Exa}[Mizohata structures] If $\nu=0$, i.e., there are no variables $z_j$, then the structure is a locally integrable Mizohata structure (see \cite{T1}), with coordinates $s,t_1,\ldots,t_{n}$ and first integral $Z(s,t)=s+i\phi(s,t)$.
\end{Exa}

A general nondegenerate structure of hypersurface type can be seen as a superposition of Mizohata structures and CR hypersurfaces of complex space.

 \subsection{The method of associated differential equations}
 The study of the relationship between the geometry of 
 real hypersurfaces in $\C^2$ and 2nd order ODEs 
  \begin{equation}\label{wzz}
 w\rq{}\rq{}=\Phi(z,w,w\rq{}).
 \end{equation}
 goes back to Segre \cite{segre}
 and Cartan \cite{cartan},
 see also Webster \cite{webster}
More generally,  the geometry of a real hypersurface in $\CC{n+1},\,n\geq 1$, is related to
 that of a complete second order system of PDEs
\begin{equation}\label{wzkzl}
w_{z_kz_l}=\Phi_{kl}(z_1,...,z_n,w,w_{z_1},...,w_{z_n}),\quad \Phi_{kl}=\Phi_{lk},\quad k,l=1,...,n,
\end{equation}
  Moreover, {\em in the real-analytic case},
  this relation becomes 
more
explicit 
by means of the Segre family.
Namely,
to any real-analytic Levi-nondegenerate hypersurface $M\subset\CC{n+1},\,n\geq 1$, one can uniquely associate a holomorphic ODE \eqref{wzz} ($n=1$) or a holomorphic PDE system \eqref{wzkzl} ($n\geq 2$),
whose solutions are precisely the Segre varieties.
The Segre family of $M$ plays a role of a 
``mediator'' between the hypersurface and the associated differential equations.  
For recent work on associated differential equations in the degenerate setting, see e.g. the papers  
  \cite{divergence, nonminimalODE} of the first author with Lamel and Shafikov.

  The associated differential equation procedure is particularly simple in the case of a Levi-nondegenerate hypersurface in $\CC{2}$. In this case the Segre family is an
 anti-holomorphic 
 2-parameter family of complex holomorphic curves. 
 It then follows from the standard ODE theory that there exists a unique ODE \eqref{wzz}, for which the Segre varieties are precisely the graphs of solutions. This ODE is called \it the associated ODE. \rm

In general, both right hand sides in \eqref{wzz},\eqref{wzkzl} appear as functions determining the $2$-jet of a Segre variety 
at a given point
as a function of the $1$-jet at the same point.
More explicitly, we use coordinates
$$
	(z,w)=(z_1, \ldots, z_n,w) \in \CC{n}\times \C = \CC{n+1}.
$$ 
Fix a 
 real-analytic
hypersurface
 $M\subset\CC{n+1}$
 passing through the origin, 
 and choose a 
sufficiently small neighborhood $U$
 of the origin.
  When $M$ is Levi-nodegenerate, we can associate a 
{\em complete second order system of holomorphic PDEs}
 to $M$,
which is uniquely determined by the condition that the differential equations are satisfied by all the
graphing functions $h(z,\zeta) = w(z)$ of the
Segre varities $\{Q_\zeta\}_{\zeta\in U}$ of $M$ in a
neighbourhood of the origin.

To be more precise, we consider a
so-called {\em  complex defining
 equation } (see, e.g., \cite{ber}),
$$
w=\rho(z,\bar z,\bar w),
$$ 
of $M$ near the origin, which one
obtains by substituting 
$u=\frac{1}{2}(w+\bar w),\,v=\frac{1}{2i}(w-\bar w)$ 
into 
a real-analytic defining equation and
solving for $w$ as function of $(z,\bar z, \bar w)$
by the implicit function theorem.
 The Segre
variety $Q_x$ of a point 
$$x=(a,b)\in U,\quad
a\in\CC{n},\,b\in\CC{},$$ 
is  now given
as the graph of the function
\begin{equation} 
\label{segredf}
w (z)=\rho(z,\bar a,\bar b), 
\end{equation}
where we slightly abuse the notation
using the same letter $w$ for 
both the last coordinate and
the function $w(z)$ defining a Segre variety.
Differentiating \eqref{segredf} we obtain
\begin{equation}\label{segreder} 
	w_{z_j}=\rho_{z_j}(z,\bar a,\bar b),
	\quad
	j=1,\ldots,n. 
\end{equation}
Considering \eqref{segredf} and \eqref{segreder}  as a holomorphic
system of equations with the unknowns $\bar a,\bar b$, 
in view of the Levi-nondegeneracy of $M$,
an
application of the implicit function theorem yields holomorphic functions
 $A_1,...,A_n, B$ such that
 \eqref{segredf} and \eqref{segreder} are solved by
$$
	\bar a_j=A_j(z,w,w'),\quad
	\bar b=B(z,w,w'),
$$
where we write
$$
	w' = (w_{z_1},  \ldots, w_{z_n}).
$$
The implicit function theorem applies here because the
Jacobian of the system coincides with the Levi determinant of $M$
for $(z,w)\in M$ (\cite{ber}). Differentiating \eqref{segredf} twice
and substituting the above solution for $\bar a,\bar b$ finally
yields
\begin{equation}\label{segreder2}
w_{z_kz_l}=\rho_{z_kz_l}(z,A(z,w,w'),
B(z,w,w'))=:\Phi_{kl}(z,w,w'),
\quad
k,l=1, \ldots, n,
\end{equation}
or, more invariantly,
i.e.\ independent of the coordinate choice,
\begin{equation}\label{segreder2'}
	j^2_{(z,w)} Q_x = \Phi(x, j^1_{(z,w)} Q_x).
\end{equation}
Now \eqref{segreder2}
(or \eqref{segreder2'})
is the desired complete system of holomorphic second order PDEs
denoted by $\mathcal E = \mathcal{E}(M)$.

 \begin{Def}\label{PDEdef}
 We call $\mathcal E = \mathcal{E}(M)$  \it the system of PDEs 
 associated with $M$. \rm  
\end{Def}

\subsection{Condition (E)}
Our main results rely on a 
new
{\em analytic regularizability} condition
for hypersurfaces of low regularity.
In the $C^\infty$ case,
such condition was
 introduced in \cite{koss_zaitsev} and  called    
  {\em Condition (E) (``E'' for extension)}.
 We shall give below both an invariant 
 and a coordinate-based formulations of it. 
 For the basic concepts in CR Geometry (such as Segre varieties and formal submanifolds) we refer to \cite{ber}, and for {\em jet bundles} and related concepts to \cite{CS}. 

Let 
$$
	\pi \colon J^{1,n}
	\to \CC{n+1},
$$ 
be the bundle
of $1$-jets of complex hypersurfaces of $\CC{n+1}$, 
which is a projective holomorphic 
 fiber bundle over $\CC{n+1}$ with the fiber 
isomorphic to $\mathbb P^n$,
and let
$M\subset\CC{n+1},\,n\geq 1$, 
be (for the moment) a ($C^\infty$) smooth
{\em  strictly pseudoconvex} real hypersurface. 
Then the complex tangent bundle $T^{\CC{}}{}M$ induces
the natural smooth (global) embedding 
$$
	\varphi\colon M\to J^{1,n}, \quad x \mapsto \bigl(x,T^{\CC{}}_xM]\bigr).
$$ 
The image 
$$
	\varphi(M)=:M_J \subset J^{1,n}
$$ 
 is a  smooth totally real
 $(2n+1)$-dimensional real submanifold  in the $(2n+1)$-dimensional complex manifold $J^{1,n}$
 by an observation of
 Webster \cite{webster}.
Next, 
associated with $M$ is the  smooth (weakly) pseudoconvex real hypersurface
$$
	\widehat M:=\pi^{-1}(M) \subset J^{1,n}.
$$
The manifold $M_J$ is a smooth real submanifold in $\widehat M$.  
Note that  $\widehat M$ itself is locally CR-equivalent 
to $M\times \CC{n}$ (and thus is {\em holomorphically degenerate}, see \cite{ber}).  
In what follows we denote by $U^+$ 
the {\em pseudoconvex side of $M$},
by which we mean
the subset of an open neighborhood
$U$
of a point of $M$ 
defined by $\rho<0$, i.e.\
$$
U^+ = U \cap \{\rho<0\},
$$
where $\rho$ is a local defining function of $M$ in $U$
with $d\rho\ne0$ and the
 complex hessian satisfying
$\partial\bar \partial\rho(X,\bar X)\ge0$
for $X\in T^{10}M$.
Given such $U^+$, 
we write
 $$\widehat U^+:=\pi^{-1}(U^+).$$ 
Then $\widehat M$
is obviously (weakly) pseudoconvex
and $\widehat U^+$ is the pseudoconvex side of $\widehat M$.
Since all our considerations are local,
the exact choice of neighborhoods
of the reference point won't play any role.

We next fix a point $p\in M$
along with the corresponding point
$$
\widehat p:= (\pi|_{M_J})^{-1}(p) \in M_J \subset \widehat M\subset J^{1,n}.
$$
Since $M$ is smooth, we may consider at each point $q\in M$ 
near $p$, its formal complexfication at $q$ as a formal complex hypersurface in $\CC{n+1}\times\overline{\CC{n+1}}$ obtained by complexifying the formal Taylor series
at $q$ of the defining function $\rho$.
In this way, the formal Segre variety $Q_q$ of $M$ 
is defined as a formal complex hypersurface
defined by power series in $\CC{}[[Z-q]]$,
where $Z=(z,w)\in \CC{n}\times\C$ as before.
Then the map
$2$-jets 
\begin{equation}\label{2jets}
q\in M \mapsto j^2_qQ_q \in J^{2,n}
\end{equation} 
induces a {\em smooth embedding}
of $M$ 
 into the bundle  
$$J^{2,n}=J^{2,n}(\CC{n+1})$$ 
of $2$-jets of complex hypersurfaces in $\CC{n+1}$. The space $J^{2,n}$  is canonically a fiber bundle  
$$
\pi^2_1 \colon J^{2,n}\to J^{1,n}
$$
over the $1$-jet bundle $J^{1,n}$.
The $2$-jet embedding
\eqref{2jets}
 defines a canonical section of $\pi^2_1$, 
\begin{equation}
\label{s12}
s_1^2\colon M_J\to J^{2,n}.
\end{equation}

Now we recall
 our analytic regularizability condition for a  smooth strictly pseudoconvex hypersurface \cite{koss_zaitsev}:

\begin{Def}\label{coorfree} 
We say that $M$ {\em satisfies Condition (E) at $p\in M$}
 if
 the canonical section $s_1^2$
 given by \eqref{s12} extends 
 holomorphically and smooth up to the boundary
 to a neighborhood of $\\widehat p$
 in
  the pseudoconvex side 
 $\widehat U^+ \cup \widehat M$.
\end{Def}


We next give an (equivalent to the above) coodinate-based formulation of Condition (E). 
Let $M\subset\CC{n+1}$ be a smooth hypersurface with the defining equation 
\begin{equation}\label{rhoeq}
\rho(Z,\bar Z)=0,\quad Z=(z,w)=(z_1,...,z_n,w)\in\CC{n+1},
\end{equation}
and
$p\in M$ a fixed point where
$\rho_w(p,\bar p)\neq 0$. 
Then the formal Segre variety at a point $q=(\tilde q,q_{n+1})\in M$ near $p$ is a graph of a function $w(z)$ (considered as a formal power series in $(z-\tilde q)$). Then the $2$-jets \eqref{2jets} are represented by either the scalar function 
$\Phi$ defined pointwise as $w''(z)$ for $z=\tilde q$  (case $n=1$),  or the symmetric matrix function $\Phi=(\Phi_{ij}),\,\,i,j=1,...,n$, defined pointwise as the collection of the partial derivatives
 $w_{z_iz_j}$ for $z=\tilde q$ (the general case $n\ge1$). It is possible to verify that, in turn, for $n=1$ we have
\begin{equation}\label{Phi1}
	\Phi=\frac{1}{(\rho_w)^3}
	\begin{vmatrix} \rho & \rho_z & \rho_w\\ \rho_z & \rho_{zz} 
	& \rho_{zw} \\ \rho_w & \rho_{zw} & \rho_{ww} \end{vmatrix},
\end{equation}
and for $n>1$ we have 
\begin{equation}\label{Phi}
	\Phi_{ij}=\frac{1}{(\rho_w)^3}
\begin{vmatrix} \rho & \rho_{z_j} & \rho_w\\ \rho_{z_i} & \rho_{z_iz_j} & \rho_{z_iw} \\ \rho_w & \rho_{z_jw} & \rho_{ww} \end{vmatrix}, \quad i,j=1,...,n.
\end{equation}
 (To obtain \eqref{Phi1},\eqref{Phi}, one has to differentiate the identity \eqref{rhoeq} once, assuming $w$ to be a function of $z$, express all the  $w_{z_j}$ in terms of the $1$-jet of $\rho$,
 and then differentiate \eqref{rhoeq} once more to obtain $w_{z_iz_j}=\Phi_{ij}$ in terms of the $2$-jet of $\rho$.)
  Both the scalar function \eqref{Phi1} and the matrix valued function \eqref{Phi} can be considered as either smooth functions on the strictly pseudoconvex hypersurface $M$ or as that on the totally real submanifold $M_J\subset J^{1,n}$ introduced above.

In terms of the $\Phi$-function \eqref{Phi1}-\eqref{Phi}, 
Condition (E) reads in an equivalent form as follows:
\begin{Def}\label{coor} 
We say that $M$ {\em satisfies Condition (E) at $p$}, if 
 the function $\Phi$ defined on $M_J$ by either \eqref{Phi1} or \eqref{Phi} extends  
   holomorphically and smoothly up to the boundary
   to a neighborhood of $\widehat p$ in
 the pseudoconvex side $\widehat U^+ \cup \widehat M$.
\end{Def}
It is obvious that \autoref{coor} is equivalent to \autoref{coorfree}. Now the analytic regularizability of a hypersurface can be characterized in terms of Condition (E), as stated in Theorem \autoref{koss_zaitsev_thm}.

\section{Special coordinates for non-degenerate structures}\label{sec3}

Let $(M,\mathcal{V})$ be a non-degenerate locally integrable structure of hypersurface type. We denote by $\Sigma^\circ$ the closure of $T^{\circ}\setminus \{0\}$ in $T^{\ast}M$. This is a conic subset, and we denote by $\Sigma$ the image of $\Sigma^{\circ}$ by the standard projection $\pi:T^{\ast}M\to M$. This is precisely the set of non-elliptic points of $\mathcal{V}$. We shall determine equations for these points in terms of the local coordinates of the previous section.

\begin{Thm}\label{thm31} Let $(M,\mathcal{V})$ be a smooth, non-degenerate structure of hypersurface type. Then, $\Sigma^\circ \subset T^{\ast}M$ is a conic symplectic submanifold of dimension $2\nu+2$, whose base projection $\Sigma$ is a closed $(2\nu+1)$-dimensional submanifold of $M$. In the local coordinates \ref{local_coord}, the set $U\cap \Sigma$ is described by the equations 
\[
\frac{\partial \phi}{\partial t_j}(z,\brr{z},s,t)=0,\,\,\,j=1,\ldots,n'.
\]
\end{Thm}
\begin{proof} We shall first prove that $\Sigma$ is a smooth $(2\nu+1)$-dimensional manifold. We work in the coordinates \ref{local_coord}. Assume that $(z,\brr{z},s,t) \in \Sigma \cap U$. Then, there is a non-trivial linear combination
\[
0\not= v=\sum_{j=1}^\nu (a_j+ib_j)\dd z_j + (\mu+i\nu)\left(\dd s + i \dd \phi\right) \in T'_{(z,\brr{z},s,t)},\,\,a_j,b_j,\mu,\nu \in \R,
\]
such that $\IM v = 0$. Computing this imaginary part, we obtain
\[
\IM v = \sum_{j=1}^\nu (b_j + \mu \phi_{x_j})\dd x_j + \sum_{j=1}^\nu (a_j+\mu\phi_{y_j})\dd y_j + \sum_{l=1}^{n'}\mu \phi_{t_{l}} \dd t_l + (\nu+\mu\phi_s)\dd s.
\]
Since $v\not=0$, we must have have $\mu \not=0$, which entails $\phi_{t_l}=0$ for all $l=1,\ldots,n'$. We conclude that $\dd_t \phi(z,\brr{z},s,t)=0$. Conversely, suppose that $\dd_t \phi(z,\brr{z},s,t)=0$. Then, one can take $\mu=1$, $\nu=-\phi_s(z,\brr{z},s,t)$, $a_j = -\phi_{y_j}(z,\brr{z},s,t)$ and $b_j = -\phi_{x_j}(z,\brr{z},s,t)$. Therefore, the covector
\[
v=-2i\sum_{j=1}^\nu\phi_{z_j}(z,\brr{z},s,t) \dd z_j + (1-i\phi_s(z,\brr{z},s,t)) \dd w
\]
is non-trivial and real, so $(z,\brr{z},s,t)\in \Sigma\cap U$.
\vspp

Since $\mathcal{V}$ is non-degenerate, the matrix $\left(\frac{\partial^2\phi}{\partial t_l \partial t_{l'}} \right)$, $1\leq l,l' \leq n'$, is non-singular, which implies the map
\[
(z,\brr{z},s,t) \mapsto \left(z,\brr{z},s,\partial_t \phi \right)
\]
is a diffeomorphism of a neighborhood of $0$ in $U$ onto another such neighborhood. We conclude that $\Sigma$ is a smooth submanifold of dimension $N-n' = 2\nu+n'+1-n' = 2\nu+1$. It is also closed, since the property of being elliptic is an open condition. 

\vspp

From non-degeneracy, the forms 
\[
\dd x_j,\dd y_j,\dd s,\,\,\dd(\partial_{t_l}\phi),\,\,j=1,\ldots,\nu,\,\,l=1,\ldots,n'
\]
make up a basis of $T^{\ast}M$ near an open neighborhood $U$ of $p\in M$. If $\sigma$ denotes the tautological one-form of $T^\ast M$, we can write 
\[
\sigma = \sum_{j=1}^\nu a_j \dd x_j + b_j \dd y_j + \xi \dd s + \dd \left(\partial_t \phi \cdot \eta \right).
\]
The equations for $\Sigma^0$ in $\pi^{\ast}U$ are given by $\partial_t \phi = \eta = 0$, so the pullback of $\sigma$ to $\Sigma^0$ is given by $\sum_{j=1}^\nu a_j \dd x_j + b_j \dd y_j + \xi \dd s$. Since the $2$-form
\[
\sum_{j=1}^\nu \dd a_j \wedge \dd x_j + \dd b_j \wedge \dd y_j + \dd \xi \wedge \dd s
\] 
is non-degenerate, $\Sigma^0$ is a symplectic submanifold of $T^{\ast}M$ (of dimension $2\nu+2$).
\end{proof}

The submanifold $\Sigma$ inherits an involutive structure from $\mathcal{V}$. We briefly recall the main notions (see Sections I.3 and I.4 in \cite{T1}):
\begin{Def} Let $(M,\mathcal{V})$ be a smooth locally integrable structure and $\Sigma \subset M$ an embedded smooth submanifold. Assume that the dimension of 
\[
\C T_p \Sigma \cap \mathcal{V}_p
\]
is constant for $p \in \Sigma$. Then, we say that $\Sigma$ is \textit{compatible} with the structure $\mathcal{V}$, and
\[
\mathcal{V}_\Sigma = \bigcup_{p\in M}\left(\C T_p \Sigma \cap \mathcal{V}_p\right)
\]
defines a smooth locally integrable structure $\Sigma$, which we call the \textit{induced locally integrable structure}.
\end{Def}
\begin{Rem} If $\Sigma$ is a compatible submanifold, then the orthogonal $$T'_{\Sigma} = \bigcup_{p\in M}T'_{\Sigma,p} \subset \C T^\ast \Sigma$$ is given by 
\[
T'_{\Sigma,p}=\left\{\omega\big|_{\C T_p \Sigma};\,\omega \in T'_p \right\}.
\]
\end{Rem}
A related notion is that of a \textit{non-characteristic submanifold}:
\begin{Def} Let $(M,\mathcal{V})$ be a smooth involutive structure and $\Sigma \subset M$ an embedded smooth submanifold. We say $\Sigma$ is \textit{non-characteristic} if, for every $p\in \Sigma$, the map
\begin{align*}
T^{\circ}_p &\rightarrow T^{\ast}_p \Sigma \\
v&\mapsto v\big|_{T_p \Sigma}
\end{align*}
is injective.
\end{Def}
We now apply these ideas to the structure $\mathcal{V}$ at hand. Since the vector fields $\partial/\partial t_j$ are transversal to the submanifold $\Sigma$, we can write it locally as graph $t=F(z,\brr{z},s)$, for some smooth function $F$ fixing the origin. Introducing the change of variables $$t\mapsto t-F(z,\brr{z},s),$$ the submanifold $\Sigma$ is defined, near $0$, by the equations $$t_j=0,\,\,\, j=1,\ldots,n'.$$ In these new coordinates, a system of first integrals for $\mathcal{V}$ is given by
\begin{equation}\label{local_coord2}
z_j=x_j+iy_j,\,\,\,j=1,\ldots,\nu,\,\,\,\,\,\,\,\,\,\,w=s+i\phi(z,\brr{z},s,t),
\end{equation}
where $\phi(0)=0$, $\dd \phi(0)=0$ and $\phi_{t_j}(z,\brr{z},s,0)\equiv 0$ for all $j=1,\ldots,n'$ and sufficiently small $(z,\brr{z},s)$.

\vspp

Now we obtain the following:
\begin{Prop} Let $(M,\mathcal{V})$ be a smooth non-degenerate structure of hypersurface type. Let $\Sigma$ be the submanifold constructed in \ref{thm31}. Then, $\Sigma$ is a non-characteristic compatible submanifold of $M$. Moreover, the induced involutive structure on $\Sigma$ is that of a locally integrable Levi-nondegenerate CR structure of hypersurface type.
\end{Prop}
\begin{proof} Let $p\in \Sigma$. Working in the coordinates \ref{local_coord2} centered at $p$, we have $\dd s$ as a basis for $T^{\circ}_p$, and clearly the curve $s\mapsto (0,0,s,0)$ is contained in $\Sigma$. In particular, the restriction $\dd s\big|_{T_p\Sigma}$ does not vanish identically, which implies $\Sigma$ is non-characteristic. Compatibility follows from the expression of the vector fields \ref{vector_fields} (which remain of the same form after the coordinate change \ref{local_coord2}), from which we see that $L_{j}$, $j=1,\ldots,\nu$, are precisely (complex) tangent to $\Sigma$. In particular, the structure is CR Levi-nondegenerate of hypersurface type.

\end{proof}

Applying Morse's lemma (with parameters), one can simplify the coordinates \ref{local_coord2} even further.

\begin{Thm}\label{thm32} Let $(M,\mathcal{V})$ be a non-degenerate smooth structure of hypersurface type. Then, given a point $p\in M$, there are coordinates $(z,\brr{z},s,t)$ centered at $p$, a non-degenerate quadratic form $\mathcal{L}$ in $\R^{n'}$ and smooth functions 
\begin{equation}\label{goodcoord}
z_j = x_j+iy_j,\,\,\,j=1,\ldots,\nu,\,\,\,\,\,\,\,w=s+i\left(\phi(z,\brr{z},s) + \mathcal{L}(t) \right)
\end{equation}
with $\phi(0)=0$, $\dd \phi(0)$ such that $\left\{\dd z_1,\ldots,\dd z_\nu, \dd w\right\}$ span $T'$ in a neighborhood of $0$. Moreover, the structure $\mathcal{V}$ induces on the submanifold $\Sigma$ given by $\{t=0\}$ an invariantly defined locally integrable Levi-nondegenerate CR structure of hypersurface type, equivalent to the hypersurface given by the defining function
\[
\rho(z,\brr{z},w,\brr{w}) = \IM w - \phi(z,\brr{z},\RE w)
\]
in $\C^{\nu+1}$.
\end{Thm}
\begin{proof} Again, we work locally as in \ref{local_coord2}. We can apply Morse's lemma with parameters (see Lemma C.6.1 in\cite{Horm3}). If the neighborhood of $0$ is sufficiently small, the only solution of $\dd_t \phi (z,\brr{z},s,t)=0$ is $t=0$. Therefore, there is a smooth family of diffeomorphisms $t\mapsto G(z,\brr{z},s,t)$ (with $(z,\brr{z},s)$ as parameters) and a non-degenerate quadratic form $\mathcal{L}$ in $\R^{n'}$ such that
\[
\phi(z,\brr{z},s,t) = \phi(z,\brr{z},s,0) + \mathcal{L}(G(z,\brr{z},s,t)).
\]
Considering the change of variables $$(z,\brr{z},s,t)\mapsto (z,\brr{z},s,G(z,\brr{z},s,t)),$$ we obtain the result.
\end{proof}

The CR submanifold $\Sigma$ will play an important role in what follows, so we give it a special name.
\begin{Def} Let $(M,\mathcal{V})$ be a nondegenerate locally integrable structure of hypersurface type. We call $\Sigma$ the \textit{central CR manifold} of $\mathcal{V}$. Given $p \in \Sigma$, we define $m(t)$ to be the signature of the quadratic form $\mathcal{L}(t)$.
\end{Def}

\section{The central manifold and the regularizability}\label{necessary_section}

The notion of mappings between integrable structures we consider in this work is the following:
\begin{Def} Let $(M,\mathcal{V})$ and $(N,\mathcal{W})$ be two formally integrable structures over smooth manifolds $M,N$. A \textit{mapping} (or morphism) between them is a smooth map $f:M\to N$ such that $(\dd f)_p \mathcal{V}_p \subset \mathcal{W}_{f(p)}$ for all $p\in M$.
\end{Def}
\begin{Rem} We use the same notation $\dd f$ for the complexification $1\otimes \dd f$ of $\dd f$ acting on the complexified tangent spaces.
\end{Rem}
In the locally integrable case, we can restate this notion in terms of the first integrals.
\begin{Prop}\label{morphism_equiv} Let $(M,\mathcal{V})$ and $(N,\mathcal{W})$ be locally integrable structures over smooth manifolds $M,N$. Let $f:M\to N$ be a smooth map. Then, $f$ is a mapping between $\mathcal{V}$ and $\mathcal{W}$ if and only if for every solution $u\in C^\infty(U)$ of $\mathcal{W}$ (defined on an open subset $U\subset N$), the pullback $u\circ f$ is a solution of $\mathcal{V}$ over $f^{-1}(U)$.
\end{Prop}

We make the simple remark that a diffeomorphism $f:M\to N$ which is a mapping from $\mathcal{V}$ to $\mathcal{W}$ is not, in general, an equivalence of integrable structures. It is the case, however, if both of them have the same rank.

Now we specialize to the case of nondegenerate structures of hypersurface type. The first observation is that a map between structures induces a CR map between the central CR hypersurfaces.
\begin{Prop}\label{CR_map} Let $(M_1,\mathcal{V}_1)$ and $(M_2,\mathcal{V}_2)$ be two non-degenerate structures of hypersurface type, and let $\Sigma_i\subset M_i$ denote the central CR structure associated to $\mathcal{V}_i$, $i=1,2$. If $f:M_1 \to M_2$ is a smooth mapping between $\mathcal{V}_1$ and $\mathcal{V}_2$, then $f$ restricts to a CR map $\Sigma_1\to \Sigma_2$.
\end{Prop}
\begin{proof} It is clear that $f(\Sigma_1) \subset \Sigma_2$, since the derivative of $f$ is a real mapping. Moreover, since $\mathcal{V}_{\Sigma_i,p}=\C T_p \Sigma_i \cap (\mathcal{V}_i)_p$, the derivative of $f$ maps $\mathcal{V}_{\Sigma_1,p}$ into $\mathcal{V}_{\Sigma_2,f(p)}$, and is therefore a morphism from $\mathcal{V}_{\Sigma_1}$ into $\mathcal{V}_{\Sigma_2}$. Since they are CR structures, this is precisely the definition of a CR map.
\end{proof}

We can use the central CR hypersurface and the signature $m(t)$ to give a necessary condition for equivalence of non-degenerate structures of hypersurface type.
\begin{Thm}\label{nec_th0} Let $(M_i,\mathcal{V}_i)$ be a (resp. smooth or analytic) nondegenerate locally integrable structure of hypersurface type, with central CR hypersurfaces $\Sigma_i$, $i=1,2$. Let $p \in \Sigma_1$, $q\in \Sigma_2$. Then, if there is a (resp. smooth or analytic) local equivalence between $\mathcal{V}_1$ and $\mathcal{V}_2$ (from a neighborhood of $p$ to a neighborhood of $q$), then
\begin{enumerate}
\item $\Sigma_1$ is locally CR equivalent to $\Sigma_2$ (from a neighborhood of $p$ in $\Sigma_1$ to one of $p_2$ in $\Sigma_2$).
\item $m(p)=m(q)$.
\end{enumerate} 
\end{Thm}
\begin{proof} This follows immediately from Proposition \ref{CR_map}. 
\end{proof}

It turns out that the converse is also true {\em for definite structures}.
\begin{Thm}\label{suf_th0} Let $(M_i,\mathcal{V}_i)$ be a (resp. smooth or analytic) nondegenerate locally integrable structure of hypersurface type, with central CR hypersurface $\Sigma_i$, $i=1,2$. Let $p \in \Sigma_1$, $q\in \Sigma_2$. Assume furtheremore that both $(M_i,\mathcal{V}_i)$ are positive. Then, if there is a (resp. smooth or analytic) local equivalence $F$ between $(\Sigma_1,p)$ and $(\Sigma_2,q)$. then there it is the restriction onto $\Sigma_1$ of a (resp. smooth or analytic) local equivalence $\widetilde F$ between germs of $\mathcal{V}_1$ and $\mathcal{V}_2$ at $p,q$ respectively.
\end{Thm}
\begin{proof} We deal with the case of smooth $F$, and discuss the analytic one in the end of the proof. 

In view of the previous sections, we can obtain coordinates of the following kind for $\mathcal{V}_1$: there is an open set $$U = W \times I \subset \R^{2\nu+1+n'},$$ where $W\subset \R^{2\nu+1}$ is an open neighborhood of $0$ and $I=(-\delta,\delta)^{n'} \subset \R^{n'}$. According to \eqref{goodcoord}, we have coordinates $(x,y,s,t)$ in $W\times I$, $z=x+iy\in\CC{\nu}$, $s\in\mathbb R\,\,t\in I$, such that the first integrals of $\mathcal{V}_1$ are given by
\[
z=z(x,y,s,t) = x+iy,\,\,\,w=s+i\left(\phi(z,\brr{z},s)+\frac{|t|^2}{2} \right),\,\,|t|^2=t_1^2+...+t_{n'}^2,
\]
where $\phi \in C^\infty(W)$ is a smooth, real-valued function such that $\phi(0)=0$ and $\dd \phi(0)=0$. The central manifold $\Sigma_1$ is characterized by setting $t=0$ in our coordinates. In this way,  $$\Sigma_1 = \{(z,w) \in \Omega;\,\IM w = \phi(z,\brr{z},\RE w)\}$$ is a strictly pseudoconvex hypersurface in $\CC{\nu+1}$, where $\Omega \subset \C^{\nu+1}$ is an appropriate neighborhood of the origin (as follows for the positivity of definicy of the Levi form). We can also think of it as the locally integrable structure on $W$ defined by the first integrals $z$ and $w=s+i\phi(z,\brr{z},s)$. 

We perform the same procedure for $\mathcal{V}_2$, which we express with a prime $'$ on top of all the objects, and assume we are given a CR equivalence $F:W\to W'$ between $\Sigma_1$ and $\Sigma_2$. Our goal is to obtain an equivalence $F^{\bullet}:U\to U'$ between $\mathcal{V}_1$ and $\mathcal{V}_2$ (maybe after contraction of the neighborhoods about $0$).

We can assume the following: let $$\Omega_+=\{(z_1,z_2)\in \Omega;\,\IM w > \phi(z,\brr{z},\RE w))\}$$ be the pseudoconvex side of $\Sigma_1$ (analogously for the second structure, with a prime $'$). Then, there is a map $$\widetilde{F}:\Omega_+ \cup \Sigma_1 \to \Omega'_+ \cup \Sigma_2$$ which is by biholomorphic as a map $\Omega_+\to \Omega'_+$, extends smoothly to $\Omega_+\cup \Sigma_1$, and its restriction to $\Sigma_1$ is the CR equivalence $\Sigma_1 \to \Sigma_2$. We can \textit{smoothly} extend $\widetilde{F}$ to be a smooth diffeomorphism between full neighborhoods of $\Sigma_1$ and $\Sigma_2$ (biholomorphic only on the side $\Omega_+$).

Let $\widetilde{F}(z,w)=(f(z,w),g(z,w))\in\CC{\nu}\times\CC{}$. In order to simplify notation, we write
\[
z:=x+iy,\,\,\,\,\,w:=s+i\left(\phi(z,\brr{z},s) + \frac{|t|^2}{2} \right),\,\,\,\,\,\zeta:=s+i\phi(z,\brr{z},s),\,\,\,(x,y,s,t)\in U.
\]
In view of the discussion above, we have the basic identity for $\widetilde F$:
$$\IM g(z,\zeta)=\phi'(f(z,\zeta),\overline{f(z,\zeta)},\RE g(z,\zeta)).$$
Note that the function $\rho(z,\brr{z},\eta,\brr{\eta}) := \IM g(z,\eta) - \phi'(f(z,\eta),\brr{f(z,\eta)},\RE g(z,\eta))$ is a \textit{smooth defining function} for the hypersurface $\Sigma_1$ in a neighborhood of the origin. In particular, it develops into the following identity:
\begin{equation}\label{developed}
\IM g(z,\eta)-\phi'(f(z,\eta),\overline{}f(z,\eta),\RE g(z,\eta))=\lambda(z,\eta)(\IM\eta-\phi(z,\brr{z},\RE \eta)),
\end{equation}
for $(z,\eta)$ in a full neighborhood of the origin, where $\lambda$ is a smooth non-vanishing function. Note that $\lambda(0)=1$, so contracting the neighborhoods about the origin allows us to assume $\lambda>0$.
Then, we set 
\[
F^{\bullet}(x,y,s,t)=\left(X(x,y,s,t),Y(x,y,s,t),S(x,y,s,t), T(x,y,s,t) \right),
\]
where 
\[
X(x,y,s,t) = \RE f\left(z,w\right),\,\,\,Y(x,y,s,t) =  \IM f(z,w),
\]
$$S(x,y,s,t) = \RE g(z,w),\quad T(x,y,s,t)=t\sqrt{\lambda(z,w)}$$
for $(x,y,s,t)\in U$ (everywhere here $w=s+i(\phi(z,\bar z,s)+|t|^2/2)$ and $z=x+iy$). Observe that $$(z,w)=(z,\zeta+i|t|^2/2) \in \Omega_+,$$  therefore the expression for $F^{\bullet}$ is well-defined. In view of the vanishing of $\lambda -1$ on $\Sigma$, $F^{\bullet}$ is smooth diffeomorphism in a neighborhood of the origin.  Now, to show it defines a map of $\mathcal{V}_1\to \mathcal{V}_2$, we have to prove that the pullback of the first integrals of $\VV_2$ are solutions of $\VV_1$. Indeed,
\[
((F^{\bullet})^{\ast}z')(x,y,s,t) = X+iY = f(z,w),
\]
which is a holomorphic function of the first integrals $(z,w)$, so it is a solution. For $w'$, we have
\[
((F^{\bullet})^{\ast} w')(x,y,s,t) = S+i\left(\phi'(X,Y,S)+\frac{|T|^2}{2} \right) = g(z,w)
\]
(as follows from \eqref{developed}), which is again a holomorphic function of the first integrals $(z,w)$. We conclude that $F^{\bullet}$ is an equivalence between $\VV_1$ and $\VV_2$.

In the analytic case, following the proof, we easily see that $\lambda$ above can be chosen analytic (since the extension $\widetilde F$ becomes actually holomorphic by the reflection principle), thus the equivalence map $F^\bullet$ is analytic too.
\end{proof}

Theorems \autoref{nec_th0} and \autoref{suf_th0} now imply the following result, which is interesting even independently of the regularizability problem.

\begin{Thm}\label{equivalence}
Let $(M_i,\mathcal{V}_i)$ be a (resp. smooth or analytic) nondegenerate locally integrable structure of hypersurface type, with central CR hypersurfaces $\Sigma_i$, $i=1,2$. Let $p \in \Sigma_1$, $q\in \Sigma_2$. Assume furtheremore that both $(M_i,\mathcal{V}_i)$ are positive. Then, the structures $\mathcal{V}_1$ and $\mathcal{V}_2$ at $p,q$ respectively are (resp. smoothly or analytically) equivalent if and only if the germs  $(\Sigma_1,p)$ and $\Sigma_2,q)$ of their central manifolds are (resp. smoothly or analytically) equivalent
\end{Thm}

Theorem \autoref{equivalence} immediately implies Theorem \autoref{main} from the Introduction, characterizing the analytic regularizability of structures in terms of their central manifolds. 

We give below an example of using the central manifold for investigating the regularizability of nondegenerate structures.

\begin{Exa}
Consider the structure $\mathcal V$ defined in a neighborhood of $\R^4$ (with variables $(x,y,s,t)$), with the first integrals
$$Z=z, \quad W=s+i(h(|z|^2)+t^2),\quad h\in C^{\omega}(U\setminus\{0\})\cap C^\infty(U),\,\,h\not\in C^{\omega}(U),\quad U=(-a,a)\subset\mathbb R$$ where the function $h$ satisfies
$$h(0)=h'(0)=0,\quad h''(0)\neq 0.$$
In this case, $\nu=n'=1$ and the vector fields that span the structure are 
\[
L_1 = \frac{\partial}{\partial \brr{z}}-izh'(|z|^2)\frac{\partial}{\partial s},\,\,\,\,L_2 = \frac{\partial}{\partial t}-2it\frac{\partial}{\partial s}.
\]
The central manifold $\Sigma$ is given by
$$\{\IM w=h(|z|^2)\}\subset\CC{2}.$$
Following the elimination procedure in Section 2, we conclude that the Segre varieties of $\Sigma$ near points with $z\neq 0$ satisfy the ODE
$$z^2w''=\psi(zw'),$$
where the function $$\psi\in C^{\omega}(U\setminus\{0\})\cap C^\infty(U),\quad \psi(0)=\psi'(0)=0,$$ is related with $h$ via
$$\xi^2h''(\xi)=\psi(\xi h'(\xi)),\quad \xi\in U.$$
Considering the latter identity as an ODE for $h$ with $\psi$ fixed, we conclude that $\psi\not\in C^{\omega}(U)$ (otherwise $h$ would be analytic at $0$ too as a solution of an analytic ODE). This means that $\Sigma$ does not satisfy Condition (E), and so $V$ is not analytically regularizable by Theorem \autoref{main}.
\end{Exa}


\section{Additional sufficient conditions: The external CR structure}

In this section, we work on sufficient conditions to ensure analytic regularizability of $(M,\mathcal{V})$. We first describe a procedure introduced by Marson \cite{Marson} that yields, in the coordinate system just introduced, an associated (generic) CR manifold to $(M,\mathcal{V})$. The main idea is to formally add new variables $x^{\bullet}$ that will turn the non-CR variables $t_j$ into imaginary parts of new complex variables. 

We maintain the notation of the previous section. Let $X\subset \R^{n'}$ be a neighborhood of $0$ and let $U\subset \R^N$ be an open neighborhood of the origin for which we have the coordinates described in Sections \ref{preliminaries} and \ref{sec3}. Set $U^{\bullet}=X\times U \subset \R^{n'+N}$. Let $M^{\bullet}=Z^\bullet(U^\bullet)$ be the submanifold parametrized by the mapping
\[
Z^{\bullet}: U^{\bullet} \to \C^{m+n'}
\]
given by
\[
\begin{cases}
\ds Z^{\bullet}_j(x,y,x^{\bullet},t,s) = x^{\bullet}_j+it_j,\,\,j=1,\ldots,n',\\
\\
\ds Z^{\bullet}_{n'+j}(x,y,x^{\bullet},t,s)=Z_j(x,y) = x_j+iy_j,\,\,j=1,\ldots,\nu,\\
\\
\ds W^{\bullet}(x,y,s,t)=W(x,y,s,t)=s+i\phi(x,y,s,t).\,\,
\end{cases}
\]
This map is easily seen to have full rank, and induces a CR structure $\mathcal{V}^{\bullet}$ on $U^\bullet$, with dual structure bundle $T'^{\bullet}$ being generated by $$\left\{\dd Z^{\bullet}_1,\ldots,\dd Z^\bullet_{m+n'}\right\}.$$ We write $z^{\bullet}_l = x^{\bullet}_l + i t_l$ and $z_j=x_j+iy_j$ for all $1\leq l \leq n'$ and $1\leq j \leq \nu$. Then, we can consider the vector fields 
\begin{equation}\label{cr_vec}
\begin{cases}
\ds L^\bullet_{l} = \frac{\partial}{\partial \brr{z^{\bullet}_l}} - i \frac{\phi_{\overline{z^{\bullet}_l}}}{1+\phi_s}\frac{\partial}{\partial s} ,\,\,l=1,\ldots,n',\\
\\
\ds L_{j}=\frac{\partial}{\partial \brr{z_j}}-i\frac{\phi_{\overline{z_j}}}{1+\phi_s}\frac{\partial}{\partial s} ,\,\,j=1,\ldots,\nu.
\end{cases}
\end{equation}
It's clear that $L^{\bullet}_kZ^{\bullet}_l=L_j Z^{\bullet}_l=0$ for all $1\leq k \leq n'$ and $1\leq l \leq \nu$ (we maintain, for simplicity, the name $L_j$ for the second set of vector fields, but formally they are the pullback of the vector fields $L_j$ in $U$ by the projection $X\times U \to U$). Therefore, 
\[
L^{\bullet}_1,\ldots,L^{\bullet}_{n'},L_1,\ldots,L_\nu \text{ generate }\mathcal{V}^{\bullet}\text{ over }U^\bullet.
\]
The new structure obtained has the same rank as the original structure, but we remark that it is \textit{not invariant}, meaning that it depends on the particular coordinates $(x,y,s,t)$ chosen.

We can describe the image of $Z^{\bullet}$ explicitely as
\[
M^{\bullet}=\left\{ (z^{\bullet},z,w) \in \C^{n'}\times \C^{\nu}\times \C;\,\IM w = \phi(z,\brr{z},\IM z^{\bullet},\RE w)\right\}
\]
This hypersurface is called an \textit{external CR structure} associated to $(\mathcal{V},p)$. We shall denote the Levi form of $\mathcal{V}^{\bullet}$ (at the characteristic covector $\dd s\big|_{0}$) by $\mathcal{L}^{\bullet}(\zeta,\zeta^{\bullet})$, where $\zeta \in \C^{\nu}$ and $\zeta^{\bullet} \in \C^{n'}$. From Section \ref{preliminaries}, the Levi form is represented by the matrix
\[
\begin{pmatrix}
\left(\frac{\partial^2 \phi}{\partial z_{j}\partial \brr{z_k}}\right)_{1\leq j,k\leq \nu} & \left(\frac{\partial^2 \phi}{\partial z_j \partial \brr{z^{\bullet}_r}} \right)_{\substack{1 \leq j \leq \nu \\ 1\leq r \leq n'}} \\
\left(\frac{\partial^2 \phi}{\partial z^{\bullet}_l \partial \brr{z_k}} \right)_{\substack{1 \leq l \leq n' \\ 1\leq k \leq \nu}} & \left(\frac{\partial^{2}\phi}{\partial z^{\bullet}_l \partial \brr{z^{\bullet}_r}} \right)_{1\leq l,r \leq n'}, 
\end{pmatrix}
\]
which is given by 
\[
\begin{pmatrix}
\left(\frac{\partial^2 \phi}{\partial z_{j}\partial \brr{z_k}}\right)_{1\leq j,k\leq \nu} &  \frac{i}{2} \left(\frac{\partial^2 \phi}{\partial z_j \partial t_r} \right)_{\substack{1 \leq j \leq \nu \\ 1\leq r \leq n'}} \\
-\frac{i}{2}\left(\frac{\partial^2 \phi}{\partial t_l \partial \brr{z_k}} \right)_{\substack{1 \leq l \leq n' \\ 1\leq k \leq \nu}}  & \frac{1}{4}\left(\frac{\partial^{2}\phi}{\partial t_l \partial t_r} \right)_{1\leq l,r \leq n'}.
\end{pmatrix}
\]
Therefore, if $\mathcal{L}(\zeta,\varpi)$ is the Levi form of $\mathcal{V}$ at the covector $\dd s\big|_{0}$, we conclude that $$\mathcal{L}^{\bullet}(\zeta,\zeta^{\bullet}) = \mathcal{L}(\zeta,(2i)^{-1}\zeta^{\bullet}).$$ In particular, we obtain the following
\begin{Prop} If $\mathcal{V}$ is a nondegenerate structure of hypersurface type such that the Levi form is positive at $\dd s\big|_{0}$, then the external CR hypersurface is strictly pseudoconvex at the origin.
\end{Prop}


We can now proceed with the proof of Theorem \ref{sufficient_thm}. We work in the coordinates described in Sections \ref{preliminaries} and \ref{sec3}. Recall that we have, in a neighborhood of the origin, a smooth real-valued function $\phi$ vanishing at $0$, with derivatives vanishing at $0$ such that
\[
\begin{cases}
Z_j(x,y) = x_j + iy_j=z_j,\,\,j=1,\ldots,\nu,\\
W(x,y,s,t)=s+i\phi(x,y,s,t)
\end{cases}
\]
and such that $\dd Z_j$ and $\dd W$ span, in a neighborhood of $0$, the bundle $T'$. We have a smooth, strictly pseudoconvex hypersurface 
\begin{equation}\label{cr_def}
M^{\bullet}=\left\{ (z^{\bullet},z,w) \in \C^{n'}\times \C^{\nu}\times \C;\,\IM w = \phi(z,\brr{z},\IM z^{\bullet},\RE w)\right\}
\end{equation}
which we assume satisfies condition (E). From Theorem \ref{koss_zaitsev_thm}, there is a CR diffeomorphism $f:M^{\bullet}\to N$ (near the origin), where $(N^{\bullet},q)$ is a real-analytic hypersurface in $\C^{n'+\nu+1}$. 

We make use of the Abelian subalgebra $\mathfrak{a} \subset \aut(M^{\bullet},p)$ generated by the real vector fields $$\left\{\frac{\partial}{\partial z^{\bullet}_l}+\overline{\frac{\partial}{\partial z^{\bullet}_l}}\right\}, \,\,\,l=1,\ldots,n'.$$
Its presence implies that the real-analytic hypersurface $N^{\bullet}$ has an abelian subalgebra $\til{\mathfrak{a}} \subset \hol(N^{\bullet},q)$ of real dimension $n'$. Since $N^{\bullet}$ is real-analytic, there are germs $V_1,\ldots,V_{n'}$ of holomorphic vector fields in $\C^{n'+\nu+1}$ such that 
\[
\spn\{\RE V_1,\ldots, \RE V_{n'}\} = \til{a}
\]
and $\RE V_1,\ldots,\RE V_{n'}$ are tangent to $N^{\bullet}$. Applying the holomorphic version of Frobenius' theorem, we can find holomorphic coordinates $(Z^{\bullet},Z,W) \in \C^{n'+\nu+1}$ near $q$ and a real-analytic function $\psi$ (vanishing, along with its derivatives at the origin) such that $N^{\bullet}$ is given by
\[
\left\{(Z^{\bullet},Z,W) \in \C^{n'+\nu+1}:\,\,\IM W = \psi(Z,\brr{Z},\IM Z^{\bullet},\RE W) \right\}
\]
near $q$. Moreover, since $f$ maps $\mathfrak{a}$ to $\til{\mathfrak{a}}$, we conclude that $f$ is of the form 
\[
f(z^{\bullet},z,w) = (z^{\bullet}+f_1(z,w),f_2(z,w)).
\]
We now consider the real-analytic locally integrable structure $(N,\mathcal{W})$ given by the first integrals
\[
\begin{cases}
\til{Z}_j(x,y,s,t)=x_j+iy_j,\,\,j=1,\ldots,\nu,\\
\til{W}(x,y,s,t)=s+i\psi(x,y,s,t).
\end{cases}
\]
Then, the induced map $F:M\to N$, given by 
\[
F(x,y,s,t) = \left(x,y,\RE f_2(Z(x,y),W(x,y,s,t)),t+\IM f_1(Z(x,y),W(x,y,s,t)) \right),
\]
defines an isomorphism between the structures $\mathcal{V}$ and $\mathcal{W}$, for the pullback of the first integrals $\til{Z}_j$ and $\til{W}_j$ are solutions of $\mathcal{V}$. We conclude that $M$ is equivalent (near $p$) to a real-analytic structure.

\qed


\section{An application to rigid structures}

We start with a 
\begin{Def}
A locally integrable structure     \eqref{vector_fields} is called {\em rigid}, if the respective function $\phi$ satisfies $\phi_s=0$.
\end{Def}
This definition generalizes the important class of rigid CR structures in CR geometry (see, e.g., \cite{ber}). 

We first establish the following, somewhat surprising, result on the analytic regularizability of rigid CR structures.

\begin{Thm}\label{rigidCR}
Let $M\subset\CC{n+1}$ be a rigid strictly pseudoconvex CR-hypersurface
$$M=\{\im w=H(z,\bar z)\},\quad (z,w)\in U\times\CC{},$$
where $U\subset\CC{n}$ is a domain and $H\in C^\infty(U)$. Then $M$ is analytically regularizable at a point $p\in M$ if and only if $M$ is actually analytic near $p$. 
\end{Thm}
\begin{proof}
We have to prove the necessity, since the sufficiency is trivial. We assume, without loss of generality, $p=0,\,H(0)=0,\,dH(0)=0$. According to Section 2, $M$ satisfies Condition (E), which means that the associated $\Phi$-function  \eqref{Phi1} or \eqref{Phi1} extends holomorphically to a domain $\widehat \Omega^+\subset J^1(\CC{n},\CC{})$ of the kind $\Omega^+\times G$, where $G$ is a neighborhood of the origin in $\CC{n}$ and $\Omega^+$ is a one-sided pseudoconvex neighborhood of $M$ near $0$. The fact that $M$ is rigid means that it is invariant under real shifts $w\mapsto w+a,\,a\in\mathbb R,$ hence the associated system \eqref{segreder2} has its defining function $\Phi$ independent of $w$. Hence, being constant on the intersection of each line $z=const,\zeta=const$ with $\widehat \Omega^+$ (for constant chosen close to the origin), $\Phi$ obviously extends holomorphically to a (full) neighborhood of the origin. This means that the Segre family of $M$ is analytic near $0$, and so is $M$.

\end{proof}

Combining now Theorem \ref{rigidCR} with Theorem \ref{main}, we arrive to the following regularizability theorem.

\begin{Thm}\label{rigid}
Let $\VV$ be a smooth rigid nondegenerate structure with a definite Levi form. Then $\mathcal V$ is analytically regularizable at a point if and only if its central manifold is actually analytic near $p$. 
\end{Thm}

\printbibliography

@book{BCH_book,
  title={An Introduction to Involutive Structures},
  author={Berhanu, S. and Cordaro, P. D. and Hounie, J.},
  series={New Mathematical Monographs},
  year={2014},
  publisher={Cambridge University Press}
}

@book{BER,
 author = {Baouendi, M. S. and Ebenfelt, P. and Rothschild, L. P.},
 publisher = {Princeton University Press},
 title = {Real Submanifolds in Complex Space and Their Mappings (PMS-47)},
 year = {1999}
}

@article{C,
  title={Global hypoellipticity for $\dbar_b$ on certain compact three dimensional CR manifolds},
  author={Cordaro, P. D.},
  journal={Resenhas do Instituto de Matemática e Estatística da Universidade de São Paulo},
  volume={2},
  number={4},
  year={1996}
}

@article{CordTreves,
  title = {Homology and cohomology in hypo-analytic structures of the hypersurface type},
  volume = {1},
  number = {1},
  journal = {Journal of Geometric Analysis},
  author = {Cordaro,  P. and Treves,  F.},
  year = {1991},
  pages = {39–70}
}

@article{BaouendiChangTreves,
author = {Baouendi, M. S. and Chang, C. H. and Treves, F.},
title = {{Microlocal hypo-analyticity and extension of CR functions}},
volume = {18},
journal = {Journal of Differential Geometry},
number = {3},
pages = {331 -- 391},
year = {1983}
}

@book{T1,
 author = {Treves, F.},
 publisher = {Princeton University Press},
 title = {Hypo-Analytic Structures: Local Theory},
 year = {1992}
}

@Misc{Nirenberg,
 Author = {Nirenberg, L.},
 Title = {A complex {Frobenius} theorem},
 Year = {1958},
 HowPublished = {Sem. analytic functions 1, 172-189 (1958).},
}

@article {koss_zaitsev,
    AUTHOR = {Kossovskiy, I. and Zaitsev, D.},
     TITLE = {Real-analytic coordinates for smooth strictly pseudoconvex
              {CR}-structures},
   JOURNAL = {Math. Res. Lett.},
  FJOURNAL = {Mathematical Research Letters},
    VOLUME = {29},
      YEAR = {2022},
    NUMBER = {5},
     PAGES = {1461--1484}
}

@article {Marson,
    AUTHOR = {Marson, M. E.},
     TITLE = {Wedge extendability for hypo-analytic structures},
   JOURNAL = {Comm. Partial Differential Equations},
    VOLUME = {17},
      YEAR = {1992},
    NUMBER = {3-4},
     PAGES = {579--592}
}

@book{Horm3,
  title = {The Analysis of Linear Partial Differential Operators III: Pseudo-Differential Operators},
  journal = {Classics in Mathematics},
  publisher = {Springer Berlin Heidelberg},
  author = {H\"{o}rmander,  Lars},
  year = {2007}
}

@article {kkz,
    AUTHOR = {Kolar, Martin and Kossovskiy, Ilya and Zaitsev, Dmitri},
     TITLE = {Normal forms in {C}auchy-{R}iemann geometry},
 BOOKTITLE = {Analysis and geometry in several complex variables},
    SERIES = {Contemp. Math.},
    VOLUME = {681},
     PAGES = {153--177},
 PUBLISHER = {Amer. Math. Soc., Providence, RI},
      YEAR = {2017},
   MRCLASS = {32V35 (32H02 32V40)},
  MRNUMBER = {3603888},
MRREVIEWER = {Gerd Schmalz},
}

@article {divergence,
    AUTHOR = {Kossovskiy, I. and Shafikov, R.},
     TITLE = {Divergent {CR}-equivalences and meromorphic differential
              equations},
   JOURNAL = {J. Eur. Math. Soc. (JEMS)},
  FJOURNAL = {Journal of the European Mathematical Society (JEMS)},
    VOLUME = {18},
      YEAR = {2016},
    NUMBER = {12},
     PAGES = {2785--2819},
      ISSN = {1435-9855},
   MRCLASS = {32V40 (32H02)},
  MRNUMBER = {3576537},
MRREVIEWER = {Ji\v r\'\i  Lebl},
       DOI = {10.4171/JEMS/653},
       URL = {https://doi.org/10.4171/JEMS/653},
}

@article {nonminimalODE,
    AUTHOR = {Kossovskiy, I. and Shafikov, R.},
     TITLE = {Analytic differential equations and spherical real
              hypersurfaces},
   JOURNAL = {J. Differential Geom.},
  FJOURNAL = {Journal of Differential Geometry},
    VOLUME = {102},
      YEAR = {2016},
    NUMBER = {1},
     PAGES = {67--126},
      ISSN = {0022-040X},
   MRCLASS = {32V20 (34M05 35H10 35H20)},
  MRNUMBER = {3447087},
MRREVIEWER = {Mauro Nacinovich},
       URL = {http://projecteuclid.org/euclid.jdg/1452002878},
}

@article{chern,
author = {Chern, S.S. and Moser, J. K.},
title = {{Real hypersurfaces in complex manifolds}},
journal = {Acta Math.},
year = {1974},
volume = {133},
pages = {219--271}
}

@article {webster,
    AUTHOR = {Webster, S. M.},
     TITLE = {On the mapping problem for algebraic real hypersurfaces},
   JOURNAL = {Invent. Math.},
  FJOURNAL = {Inventiones Mathematicae},
    VOLUME = {43},
      YEAR = {1977},
    NUMBER = {1},
     PAGES = {53--68},
      ISSN = {0020-9910},
   MRCLASS = {32C05 (32F99)},
  MRNUMBER = {0463482},
MRREVIEWER = {Klaus Fritzsche},
       DOI = {10.1007/BF01390203},
       URL = {https://doi.org/10.1007/BF01390203},
}

@article {cartan,
    AUTHOR = {Cartan, \'Elie},
     TITLE = {Sur la g\'eom\'etrie pseudo-conforme des hypersurfaces de l'espace
              de deux variables complexes {II}},
   JOURNAL = {Ann. Scuola Norm. Sup. Pisa Cl. Sci. (2)},
  FJOURNAL = {Annali della Scuola Normale Superiore di Pisa. Classe di
              Scienze. Serie II},
    VOLUME = {1},
      YEAR = {1932},
    NUMBER = {4},
     PAGES = {333--354},
      ISSN = {0391-173X},
   MRCLASS = {DML},
  MRNUMBER = {1556687},
       URL = {http://www.numdam.org/item?id=ASNSP_1932_2_1_4_333_0},
}

@Article{segre,
    Author = {Beniamino {Segre}},
    Title = {{Questioni geometriche legate colla teoria delle funzioni di due variabili complesse.}},
    FJournal = {{Rendiconti del Seminario Matematico delle Facolt\'a di Scienze della R. Universit\'a di Roma, II. Serie}},
    Journal = {{Rend. Sem. Mat. Roma, II. Ser.}},
    Volume = {7},
    Number = {2},
    Pages = {59--107},
    Year = {1932},
    Publisher = {Universit\'a degli Studi, Seminario Matematico, Roma},
    Language = {Italian},
    Zbl = {0005.01901}
}

@Article{KZnew,
    Author = {I. Kossovskiy, D.Zaitsev},
    Title = {{Sphericity and Analyticity of a Real Hypersurface in Low Regularity.}},
    FJournal = {{Arxiv. org}},
    Journal = {{To appear}},
    Volume = {},
    Number = {},
    Pages = {},
    Year = {2025},
 
}

@book {CS,
    AUTHOR = {\v{C}ap, Andreas and Slov\'{a}k, Jan},
     TITLE = {Parabolic geometries. {I}},
    SERIES = {Mathematical Surveys and Monographs},
    VOLUME = {154},
      NOTE = {Background and general theory},
 PUBLISHER = {American Mathematical Society, Providence, RI},
      YEAR = {2009},
     PAGES = {x+628},
      ISBN = {978-0-8218-2681-2},
   MRCLASS = {53C10 (17B60 22E60 53A55 53C30 58J70)},
  MRNUMBER = {2532439},
MRREVIEWER = {Stuart Armstrong},
       DOI = {10.1090/surv/154},
       URL = {https://doi.org/10.1090/surv/154},
}

\end{document}